\documentclass[titlepage,11pt]{article}
\usepackage{latexsym}
\usepackage{amsfonts}
\usepackage{amsmath}
\usepackage{mathtools}
\usepackage{tikz}
\usepackage{comment}
\usetikzlibrary{shapes.geometric}
\newcommand\blackslug{\hbox{\hskip 1pt \vrule width 4pt height 8pt depth 1.5pt
        \hskip 1pt}}
\newcommand\bbox{\hfill \quad \blackslug \bigbreak}

\def\LL{,\ldots,}
\newcommand{\vare}{\varepsilon}

\newcommand{\cupcup}{\cup \cdots\cup}

\def\polylog{\operatorname{polylog}}

\title{Subdivisions and near-linear stable sets}
\author{
Tung Nguyen\thanks{Supported by AFOSR grant
FA9550-22-1-0234, and by NSF grant  DMS-2154169, and a Porter Ogden Jacobus Fellowship.}\\
Princeton University,\\ Princeton, NJ 08544, USA
\and
Alex Scott\thanks{Supported by EPSRC grant EP/X013642/1}\\
University of Oxford, \\
Oxford, UK
\and
Paul Seymour\thanks{Supported by AFOSR grant
FA9550-22-1-0234, and by NSF grant  DMS-2154169.}\\
Princeton University,\\ Princeton, NJ 08544, USA}

\date{May 10, 2024; revised \today}

\newtheorem{thm}{}[section]

\newcommand{\Proof}{\noindent{\bf Proof.}\ \ }

\begin{document}
\maketitle
\begin{abstract}
We prove that for every complete graph $K_t$, all graphs $G$ with no induced subgraph isomorphic to a subdivision of $K_t$ have a stable subset of 
size at least
$|G|/\operatorname{polylog}|G|$. This is close to best possible, because for $t\ge 7$, not all such graphs
$G$ have a stable set of linear size, even if $G$ is triangle-free.

\end{abstract}

\section{Introduction}

Graphs in this paper are finite, and have no loops or parallel edges. For a graph $G$, $|G|$ denotes the number of vertices of $G$,
$\chi(G)$ is its chromatic number, and $\omega(G)$ and $\alpha(G)$ denote the sizes of its largest clique and stable set 
respectively. If $G,H$ are graphs, $G$ is {\em $H$-free} if no induced subgraph of $G$ is isomorphic to $H$, and if $\mathcal{H}$
is a set of graphs, $G$ is {\em $\mathcal{H}$-free} if $G$ is $H$-free for each $H\in \mathcal{H}$. 
We say $G$ is {\em $H$-subdivision-free} if no induced subgraph of $G$ is isomorphic to a subdivision of $H$. All logarithms are 
to base two.

We will prove:
\begin{thm} \label{stable2}
Let $t\ge 3$ be an integer, and let $c_t= (2t)^{-4(t-2)}$. Then 
every $K_t$-subdivision-free graph $G$ with $|G|\ge 2$ satisfies
$$\alpha(G)\ge \frac{c_t|G|}{(\log |G|)^{3t-5}}.$$
\end{thm}
Let us put this in some context. 
First, a {\em string graph} is the intersection graph of curves in the plane. String graphs with clique number at most $k$ are 
$K_t$-subdivision-free, where $t$ depends only on $k$, and so our result extends a result of Fox and Pach~\cite{fox} (and see also~\cite{fox2}):
\begin{thm}\label{string}
For every integer $k\ge 1$, there exists $c>0$ such that every string graph $G$ with clique mumber at most $k$
has a stable set of size at least
$|G|/(\log |G|)^c$.
\end{thm}
Our result is numerically weaker: $c$ in the Fox-Pach theorem is $O(\log k)$, while for us it is $O(k)$ (applying \ref{subdivision}). 
Nevertheless, it is 
pleasing to extend a geometric result to a much broader class of graphs.

Second, the Gy\'arf\'as-Sumner conjecture~\cite{gyarfastree, sumner} says:
\begin{thm}\label{GSconj}
{\bf Conjecture: }For every forest $H$ and every integer $k\ge 1$, there exists $c\ge 1$ such that $\chi(G)\le c$ for every $\{H,K_{k+1}\}$-free graph $G$.
\end{thm}
In particular, 
if the Gy\'arf\'as-Sumner conjecture is true, then every $\{H,K_{k+1}\}$-free graph $G$
has a stable set of size at least $|G|/c$. This is still open, but in~\cite{stable1} we proved that it is ``nearly'' true, in the following sense:
\begin{thm}\label{stable1}
For every forest $H$ and integer $k$, every $\{H,K_{k+1}\}$-free graph $G$ satisfies $\alpha(G)\ge |G|^{1-o(1)}$,
and hence has chromatic number at most $|G|^{o(1)}$.
\end{thm}

One might regard this as evidence in support of the Gy\'arf\'as-Sumner conjecture, but one should not place much faith in 
such evidence, because the result of this paper provides just as much evidence to support a conjecture that turned out to be false.
Scott~\cite{scott} suggested an analogue of \ref{GSconj}, that we could replace excluding a forest with excluding all subdivisions of any given graph:
\begin{thm}\label{scottconj}
{\bf False conjecture: }For every graph $H$, and every integer $k\ge 1$, there exists $c$ such that every $K_{k+1}$-free, $H$-subdivision-free graph $G$ satisfies $\chi(G)\le c$.
\end{thm}
Scott proved that this is true if $H$ is a forest, but 
it was shown to be false in general by Pawlik, Kozik, Krawczyk, Laso\'{n},
Micek, Trotter and Walczak~\cite{sevenpoles}. Results of 
Chalopin, Esperet, Li and Ossona de Mendez~\cite{chandeliers} and Walczak~\cite{walczak} together gave a strengthening:
\begin{thm}\label{walczak}
Let $H$ be obtained from $K_4$ by subdividing once every edge in a cycle of length four. There are infinitely many $K_3$-free and 
$H$-subdivision-free graphs $G$ such that $\alpha(G)<|G|/\log\log|G|$.
\end{thm}
Every subdivision of $K_7$ contains a triangle or a subdivision of the graph $H$ of 
\ref{walczak}, so the graphs $G$ of \ref{walczak} are $K_7$-subdivision-free, as mentioned in the abstract.

However, we can prove that all $H$-subdivision-free graphs $G$ with bounded clique number 
have a stable set of size 
$|G|/\operatorname{polylog}|G|$, and so
\ref{scottconj} is ``nearly'' true in the same sense as \ref{GSconj}, despite \ref{walczak}.
We will show:
\begin{thm} \label{subdivision4again}
Let $G$ be a graph with $|G|\ge 2$, let $k\ge 1$ be an integer with $\omega(G)\le k$, and let $H$ be a graph 
such that no induced subgraph of $G$
is a subdivision of $H$.  Then
$$\alpha(G)\ge \frac{|G|}{(2|H|)^{4(k-1)}(\log |G|)^{3k-2}}.$$
\end{thm}
This immmediately implies a polylogarithmic bound $O((\log |G|)^{3k-1})$ on the chromatic number of such graphs $G$.

We remark that if we work with induced minors rather than induced subdivisions, there is a related result, proved by T. Korhonen and D. Lokshtanov~\cite{korhonen}:
\begin{thm}\label{korhonen}
If $G$ is a graph with $m$ edges, and does not contain a graph $H$ as an induced minor, then 
$G$ has a separator of cardinality $O(\sqrt{m})$. 
\end{thm}
(We omit the definitions of ``induced minor'' and ``separator'', since we will not need them any more.)
By an argument of Fox and Pach, this implies that if in addition $G$ satisfies $\omega(G)\le k$,
then $\alpha(G)\ge |G|(\log |G|)^{-O(\log k)}$. If we just exclude induced subdivisions of $H$ rather than induced minors, 
there need not be separators of small order: indeed, graphs with maximum degree three need not have small separators
but exclude subdivisions of $K_5$.

\section{A sketch of the proof}

To prove \ref{stable2}, for inductive purposes we will prove a stronger result, proving a bound on 
$\alpha(G)$ that depends not only on $t$ but also on $\omega(G)$ and on the maximum degree $\Delta(G)$ of $G$. In addition, we find 
that the same bound holds even if we only exclude subdivisions of $K_t$ in which each edge is subdivided at least twice and at 
most $(\log |G|)^2-1$ times.
Let us say a subdivision of $H$ has {\em lengths between $p$ and $q$} if each edge of $H$ is replaced by a path of length between $p$ and $q$.
So our main result is:

\begin{thm} \label{subdivisionagain}
Let $G$ be a graph. Let $d\ge 2$ be an integer with $\Delta(G)\le d$; let $k\ge 1$ be an integer with $\omega(G)\le k$; and let $t\ge 3$ be an integer
such that no induced subgraph of $G$
is a subdivision of $K_t$ with lengths between $3$ and $(\log |G|)^2$.  Then
$$\alpha(G)\ge \frac{|G|}{(2t)^{4(k-1)}(\log |G|)^{3(k-1)}\log d}.$$
\end{thm}

We give the proof in the next section, but let us sketch it first. 
There is an initial step, that we can assume $\Delta(G)\le |G|/(\polylog|G|)$, because if there is a vertex of sufficiently large 
degree we just apply induction on the clique number to its set of neighbours and win. 

The remainder of the proof 
breaks into two parts. In the first part we find $t$ disjoint ``stars'': $t$ suitable vertices $a_1\LL a_t$, that we intend to be the high-degree vertices
of a subdivision of $K_t$ that we will try to construct. They need to be nonadjacent, and for each to have a large stable set $B_i$ of 
neighbours nonadjacent to the other $a_j$'s; and we need these sets $B_i$ to be sparse to one another.
We try to construct $a_1\LL a_t$ greedily. Say we have selected $a_1\LL a_i$ with corresponding sets $B_1\LL B_i$. Let $G'$ be the 
set of vertices that nonadjacent to $a_1\LL a_i$, and sparse to each of $B_1\LL B_i$ (the sparsity we need here depends just on $t$).
Then $G'$ still contains nearly all the vertices of $G$, and if it has maximum degree at most $\Delta(G)/2$,
we win
by inductively applying the theorem to $G'$. This is why we need the $\log d$ term in \ref{subdivisionagain},
and this is the only place we use that.

If there is a vertex $a_{i+1}$ in $G'$ with degree at least $\Delta(G)/2$, we apply the theorem inductively to its set $X$ of neighbours,
and find a large stable subset $B_{i+1}$, and they form the next star. It is critical here that the clique number of $G[X]$ is less 
than that of $G$, and since the bound we are proving depends on clique number, the large stable subset we find is bigger (by a 
$\polylog|G|$ factor) than we could guarantee in a general set of the same size as $X$.

When we have found suitable $a_1\LL a_t$ and $B_1\LL B_t$, we go to the second part of the proof. Each pair of $B_i$'s 
needs to be joined by a path, and we need to choose these paths so that their union is an induced subgraph. So, suppose we have 
chosen some of the paths, and let us try to choose the next one. We need to put aside all vertices with neighbours in the paths 
that are already chosen, and route the next path through the remaining vertices. For this we must arrange that the number of 
vertices put aside is not too big. Each vertex only disqualifies at most $\Delta(G)$ other vertices, but we must arrange that 
the total number of vertices in the paths already selected is under control; and to do this, we insist that each path must have 
length at most some fixed power of $\log|G|$,
and we find $(\log |G|)^2$ works. 

There is another disqualification issue as well: some vertices in the $B_i$'s cannot be used because they have neighbours in 
already-selected paths. So to keep this number small, we will take care to route all paths through the set of vertices that are reasonably sparse to each $B_i$. 
Since we will have
to disqualify the neighbours in $B_i$ of about $t^2(\log|G|)^2$ vertices, we need ``reasonably sparse'' to mean having 
$O(|B_i|/(\log|G|)^2)$ neighbours in each $B_i$. The number of vertices disqualified for this reason turns out to be $O((\log|G|)^2\Delta(G))$, about the same as the number disqualified for the first reason.
So there is a set $S$ say, of size about $(\log|G|)^2\Delta(G)$, of disqualified vertices, and we have to route the next path between
$B_1$ and $B_2$ (say) without using any vertex in $S$.

To do this, we look at set expansion. Start with $N_1^0 = B_1$, and define $N_1^i$ to be the set of vertices that belong to or have 
a neighbour in $N_1^{i-1}$.
Suppose we could prove that $|N_1^i|\ge (1+2/\log |G|))|N_1^{i-1}|$ for each $i$; then when $i$ is about $(\log |G|)^2/2$, $N_1^i$
contains more than half the vertices of the graph. Do the same for $B_j$, and so $N_1^i$ and $N_2^i$ intersect, and we can pick out 
a short path between $B_i, B_j$ in their union. This is the plan, but to make it work there are several things to worry about. 

First, how can we prove that subsets have large expansion in this sense? 
If there is a set $X$, of size less than $|G|/2$ and
such that the set $N$ of vertices with neighbours in $X$ (including $X$ itself) is not very big, then we can apply
the inductive hypothesis to $X$ and to $G\setminus N$, find large stable sets in each, and take their union to find a really large 
stable set in $G$ and win. So this gives a mechanism to show that $N$ is large, for any set $X$. 

Second, we have to avoid
the disqualified vertices $S$; so we should redefine $N_1^i$ to be the set of vertices, not in $S$, that have a neighbour in 
$N_1^{i-1}$, and then that issue is resolved. But to prove that $|N_1^i|\ge (1+2/\log |G|))|N_1^{i-1}|$ with the new meaning of $N_1^i$
, we would need to prove that the 
set of {\em all} neighbours of $N_1^{i-1}$ in $G$ is somewhat larger, so that even when we remove those in $S$, what remains is still large.
So we need the set of neighbours of $N_1^{i-1}$ to be something like at least $|S|+(1+2/\log |G|))|N_1^{i-1}|$. This is not a problem for the larger values of $i$,
because then $|S|$ is fixed and all the other sets are getting big. It would be a problem when $i=1$, if all we knew about the $B_i$'s
was that their size was something like $\Delta(G)$, because $S$ is much larger than this and the expansion argument would not work.  
But when $i=1$ we have an advantage,
because $B_i$ is already stable, and we found it by using the induction on clique number, as explained earlier. That gain is 
enough to make the expansion argument work again; indeed, the expansion when $i=1$ is big enough to ensure that all the subsequent
sets $N^j_1$ are comfortably large and consequently they expand nicely.

\section{The main proof}
In this section we prove our main result, which we restate:
\begin{thm} \label{subdivision}
Let $G$ be a graph. Let $d\ge 2$ be an integer with $\Delta(G)\le d$; let $k\ge 1$ be an integer with $\omega(G)\le k$; and let $t\ge 3$ be an integer 
such that no induced subgraph of $G$                                
is a subdivision of $K_t$ with lengths between $3$ and $(\log |G|)^2$.  Then 
$$\alpha(G)\ge \frac{|G|}{(2t)^{4(k-1)}(\log |G|)^{3(k-1)}\log d}.$$
\end{thm}
\Proof  We assume inductively that the result holds for all proper induced subgraphs of $G$ (for all choices of $k\ge 1$ and $t\ge 3$).
If $k=1$, the result is trivial, so we may assume that $k\ge 2$.
If $\Delta(G)\le 1$, then $\alpha(G)\ge |G|/2$
and the result holds since $(2t)^{4(k-1)}\ge 6^4\ge 2$; so we may assume that $\Delta(G)\ge 2$. Hence, by reducing $d$ 
if possible, we may assume that $\Delta(G)=d$.
To save writing, let us define $T=(2t)^4$,  $L=\log |G|$, and $D=\log d$. Then we must show that 
$\alpha(G)\ge \frac{|G|}{T^{k-1}L^{3(k-1)}D}.$ 
We suppose this is false.
We claim:
\\
\\
(1) {\em $d D^3\le |G|/T$ and $d+1\ge T^{k-1}L^{3(k-1)}D$, and consequently $d\ge 10^9\ge 2^{29}$.}
\\
\\
Let $v$ be a vertex with degree $\Delta(G)=d$, and let $X$ be the set of all neighbours of $v$; so $|X|=d$, 
and $\omega(G[X])\le k-1$. From the inductive hypothesis,
$$\alpha(G[X])\ge \frac{d}{T^{k-2}D^{3(k-2)}D}$$
and therefore this is smaller than 
$$\frac{|G|}{T^{k-1}L^{3(k-1)}D},$$
because otherwise the result would hold. This proves the first inequality of (1). 
For the second, 
since $G$ has maximum degree $d$, and hence has a stable set of size at least $|G|/(d+1)$,
it follows that 
$$\frac{|G|}{d+1}<\frac{|G|}{T^{k-1}L^{3(k-1)}D},$$
and so $d+1\ge T^{k-1}L^{3(k-1)}D$. This proves the second inequality. Since $|G|\ge d$, and $t\ge 3$ and $k\ge 2$,
it follows that $T\ge 6^4$; since $L^{3(k-1)}D \ge D^4$, it follows that $d+1\ge 6^4(\log d)^4$
and therefore $d>10^9$.
This proves (1).
\\
\\
(2) {\em If $X\subseteq V(G)$ and $|X|\le |G|/D$, then $\Delta(G\setminus X)\ge d/2$.}
\\
\\
Let $G'=G\setminus X$, and let $d'=\lfloor d/2\rfloor$, and suppose that $\Delta(G')\le d/2$. 
Thus $\Delta(G')\le d'$, and $d'\ge 2$, so from the inductive hypothesis,
$$T^{k-1}\alpha(G')\ge \frac{|G'|}{(\log |G'|)^{3(k-1)} \log d'}\ge 
\frac{(1-1/D)|G|}{L^{3(k-1)}(D - 1)}=\frac{|G|}{L^{3(k-1)}D},$$
and so the result holds, a contradiction. This proves (2).

\bigskip

If $A,B\subseteq V(G)$ are disjoint, we say $A$ is {\em $x$-sparse to $B$} if each vertex in $A$ has at most $x|B|$ neighbours in $B$,
and the smallest such $x$ is the {\em sparsity} of $A$ to $B$.
Let us say a {\em star system} in $G$ is a pair of sequences $(a_1\LL a_n)$ and $(B_1\LL B_n)$ of the same length,
with the following properties:
\begin{itemize}
\item $a_1\LL a_n\in V(G)$ are distinct and pairwise nonadjacent;
\item $B_1\LL B_n$ are pairwise disjoint stable subsets of $V(G)\setminus \{a_1\LL a_n\}$; and
\item for $1\le i\le n$, $a_i$ is adjacent to each vertex in $B_i$ and has no other neighbour in $B_1\cupcup B_n$.
\end{itemize}
We call $n$ the {\em length} of the star system, and its {\em size} is $\min(|B_i|:1\le i\le n)$ (or $|G|$, if $n=0$). 
Its {\em semi-sparsity}
is the maximum over all $i,j$ with $1\le i<j\le n$
of the sparsity of $B_j$ to $B_i$
(or $0$, if $n\le 1$).
Its {\em sparsity} is the maximum over all distinct $i,j\in \{1\LL n\}$
of the sparsity of $B_i$ to $B_j$ (or $0$, if $n\le 1$).
\\
\\
(3) {\em Let $p$ be a nonnegative integer, and let $0<q\le 1$, such that $p/q\le 400T$. 
Then $G$ contains a star system of length $p$, 
size at least $\frac{d}{2T^{k-2}D^{3k-5}}$, and semi-sparsity at most $q$.}
\\
\\
We proceed by induction on $p$, 
and so we may assume that $p\ge 1$ and $G$ contains a star system $(a_1\LL a_{p-1})$,
$(B_1\LL B_{p-1})$ of length $p-1$,
size at least $\frac{d}{2T^{k-2}D^{3k-5}}$, and semi-sparsity at most  $q$.
Let $X_1$ be the set of vertices of $G$ that 
are adjacent to one of $a_1\LL a_{p-1}$. Let $X_2$ be the set of vertices $v$ of $G\setminus X_1$ such that
for some $i\in \{1\LL p-1\}$, $v$ has at least $q|B_i|$  neighbours in $B_i$.

Thus $|X_1|\le (p-1)d\le pd/q$, and $|X_2|\le (p-1)d/q\le pd/q$ (because there are at most $|B_i|d$ edges incident with
$B_i$); and so $X=X_1\cup X_2$ has cardinality at most 
$(2p/q)d$. Since $D^2\ge 800 $ by the last statement of (1), it follows (from the first statement of (1)) that  
$$(2p/q)dD \le  (p/q)dD^3/400\le   (p/q)|G|/(400T)  \le |G|.$$
Consequently $|X|\le |G|/D$, and so by (2), $\Delta(G\setminus X)\ge d/2$. Hence there exists $a_p\in V(G)\setminus X$
with degree at least $d/2$ in $G\setminus X$. Let $C$ be the set of its neighbours in $G\setminus X$; so $|C|\ge d/2$. 
Now $G[C]$ has clique number at most $k-1$, so from the inductive hypothesis, there exists a stable subset
$B_p\subseteq C$ with 
$$|B_p|\ge \frac{|C|}{T^{k-2}(\log |C|)^{3(k-2)}D}\ge \frac{d}{2T^{k-2}D^{3k-5}}.$$
But then $(a_1\LL a_{p-1}, a_p)$ and
$(B_1\LL B_{p-1}, B_p)$ form the desired star system.
This proves (3).
\\
\\
(4) {\em Let $p$ be a nonnegative integer, and let $0<q\le 1$, such that $p^2/q\le 200T$. 
Then $G$ contains a star system of length $p$, 
size at least $\frac{d}{4T^{k-2}D^{3k-5}}$, and sparsity at most $q$.}
\\
\\
Let $q'=q/(2p)$; then by (3), $G$ contains a star system $(a_1\LL a_p)$,
$(B_1\LL B_p)$ of length $p$, 
size at least $\frac{d}{2T^{k-2}D^{3k-5}}$, and semi-sparsity at most $q'$. 
Inductively, for $i = p, p-1\LL 1$ in turn, we define $C_i\subseteq B_i$ as follows. For $i+1\le j\le p$, since 
$C_j$ is $q'$-sparse to $B_i$, at most $|B_i|/(2p)$ vertices in $B_i$ have more than $2pq'|C_j|$ neighbours in $C_j$. Hence
there exists $C_i\subseteq B_i$ with $|C_i|\ge |B_i|/2$ that is $2pq'$-sparse and hence $q$-sparse to each of $C_{i+1}\LL C_p$. Moreover, for $i<j\le p$,
$C_j$ is $2q'$-sparse and hence $q$-sparse to $C_i$, since $|C_i|\ge |B_i|/2$. This completes the inductive definition of $C_p\LL C_1$. 
Then $(a_1\LL a_{p-1}, a_p)$,
$(C_1\LL C_{p-1}, C_p)$ satisfies (4).
\bigskip
Now we show a local expansion property:
\\
\\
(5) {\em If $A\subseteq V(G)$ is stable, then there are at least $T^{k-1}L^{3(k-1)}D|A|$ vertices that belong to or have a neighbour in $A$.}
\\
\\
Let $S$ be the set of all vertices that belong to or have a neighbour in $A$. From the inductive hypothesis applied to $G\setminus S$,
we obtain a stable set, and its union with $A$ is stable; so
$$\frac{|G|-|S|}{T^{k-1}(\log(|G|-|S|))^{3(k-1)}D}+|A|<\frac{|G|}{T^{k-1}L^{3(k-1)}D}.$$
Since $\log(|G|-|S|)\le L$, this implies
$$\frac{|G|-|S|}{T^{k-1}L^{3(k-1)}D}+|A|<\frac{|G|}{T^{k-1}L^{3(k-1)}D},$$
which simplifies to the claim. This proves (5).

\bigskip

From (4), taking $p=t$ and $q=1/(4t^2)$, we deduce that, since $4t^4\le 3200t^4=200T$, $G$ contains a star system of length $t$, 
size at least $\frac{d}{4T^{k-2}D^{3k-5}}$, and sparsity at most $1/(4t^2)$. Let $(a_1\LL a_t)$,
$(B_1\LL B_t)$ be such a star system. Thus, from the second statement of (1), each $B_i$ has cardinality at least 
$$\frac{d}{4T^{k-2}D^{3k-5}}\ge \frac{dTD^3}{4T^{k-1}L^{3(k-1)}D}\ge \frac{dTD^3}{d+1}\ge 8t^4D^3.$$

Let $X_1$ be the set of vertices adjacent to one of $a_1\LL a_t$ (so $|X_1|\le td$), and let 
$X_2$ be the set of vertices $v\in V(G)\setminus X_1$ such that for some $i\in \{1\LL t\}$, $v$ has 
at least $|B_i|/(2t^2L^2)$ neighbours in $B_i$. Thus (counting edges incident with $V(G)\setminus X_1$), we have $|X_2|\le 2t^3L^2d$. 
Let $Y\subseteq V(G)$ be some subset of $V(G)$
such that:
\begin{itemize}
\item $Y\cap X_1\subseteq B_1\cupcup B_t$ and $Y\cap X_2=\emptyset$;
\item $|Y\cap X_1|\le t^2$, and $|Y|\le \frac12 t^2L^2$.
\end{itemize}
Let us call such a set $Y$ an {\em obstruction set}.
For the moment, let us fix an obstruction set $Y$, and let $X_3$ be the set of vertices of $G$ that are equal or adjacent to a vertex in $Y$. 
\\
\\
(6) {\em $|X_3|\le \frac12 t^2L^2d$, and for $1\le i\le t$, $|X_3\cap B_i|\le |B_i|/2$.}
\\
\\
The first statement is clear, since $\Delta(G)\le d$. The vertices of $X_3\cap B_i$ are of three types: those in $Y\cap B_i$, those
with a neighbour in $(Y\setminus B_i)\cap (B_1\cupcup B_t)$, and those with a neighbour in $Y\setminus (B_1\cupcup B_t)$.
Let $|Y\cap B_i|=n$; then the number of the second type is at most $(t^2-n)|B_i|/(4t^2)$ (since $B_i$ is stable and 
the star system is $1/(4t^2)$-sparse). So the number of vertices of the first or second type is at most 
$n+(t^2-n)|B_i|/(4t^2)\le |B_i|/4$ (since $|B_i|\ge 8t^4D^3\ge 4t^2$).
Since $Y$ is disjoint from $X_2$, there are at most $|Y|\cdot |B_i|/(2t^2L^2)\le |B_i|/4$
vertices of the third type. Consequently $|X_3\cap B_i|\le |B_i|/2$. This proves (6).

\bigskip

For $1\le i\le t$ and each integer $r\ge 0$, let $N^r_i$ be the union of the vertex sets of all paths of $G$ with length at most $r$
that have one end in $B_i\setminus X_3$ and have  no other vertices in $X_1\cup X_2\cup X_3$.
\\
\\
(7) {\em For $1\le i\le t$, 
$|N_i^1|\ge 3t^3L^3d$.}
\\
\\
Let $S=N_i^1\cup X_1\cup X_2\cup X_3$. Since 
$$|B_i\setminus X_3|\ge |B_i|/2\ge \frac{d}{8T^{k-2}D^{3k-5}}$$ by (6), 
and since $S$ contains all vertices that belong to or have a neighbour in $B_i\setminus X_3$,
(5) implies that 
$$|S|\ge T^{k-1}L^{3(k-1)}D|B_i\setminus X_3|\ge \frac{dT^{k-1}L^{3(k-1)}D}{8T^{k-2}D^{3k-5}}=\frac{dTL^{3(k-1)}}{8D^{3k-6}}\ge 2t^4L^3d.$$
Since 
$$|X_1\cup X_2\cup X_3|\le td +  2t^3L^2d + \frac12 t^2L^2d\le 3t^3L^2d,$$
it follows that 
$$|N_i^1|\ge 2t^4L^3d - 3t^3L^2d\ge 3t^3L^3d.$$
This proves (7).
\\
\\
(8) {\em For $1\le i\le t$ and each integer $r\ge 2$, if $|N_i^{r-1}|\le |G|/2$ then 
$$|N^r_i|\ge \left(1+2/L\right)|N_i^{r-1}|.$$
}

\noindent
Let $|N_i^{r-1}|=m$ and let $S=N_i^{r}\cup X_1\cup X_2\cup X_3$. Let $A\subseteq N_i^{r-1}$ be stable with size at least 
$\frac{m}{T^{k-1}(\log m)^{3(k-1)}D}$. Since $S$ contains all vertices that belong to or have a neighbour in $A$, (5) implies that
$$|S|\ge T^{k-1}L^{3(k-1)}D|A|\ge \frac{L^{3(k-1)}}{(\log m)^{3(k-1)}}m.$$
Since $\log m\le L-1$, we deduce that 
$$|S|\ge \left(1+1/L\right)^{3(k-1)}m\ge (1+3/L)m.$$
As before, $|X_1\cup X_2\cup X_3|\le 3t^3L^2d\le m/L$, since
$m\ge 3t^3L^3d$ (because $N^1_i\subseteq N^{r-1}_i$ and $|N^1_i|\ge 3t^3L^3d$ by (7)).
Since
$N_i^r= S\setminus (X_1\cup X_2\cup X_3)$, it follows that 
$$|N_i^r|\ge \left(1+3/L\right)m-m/L= \left(1+2/L\right)m.$$
This proves (8).

\bigskip
We recall that earlier we fixed some obstruction set $Y$, and used it to define $X_3$. We deduce that:
\\
\\
(9) {\em For every obstruction set $Y$, and for all distinct $i,j\in \{1\LL t\}$, there is an induced path of length  
at most $L^2-2$
between $B_i,B_j$ with no vertices in $X_2$ and with no internal vertices in $X_1$, such that none of its vertices belong to or have a neighbour in $Y$.}
\\
\\
Define $X_3$ as before. 
Since $\left(1+2/L\right)^{L/2}\ge 2$, it follows that 
$\left(1+2/L\right)^{L^2/2}\ge 2^L=|G|$. Let $r=\lfloor L^2/2\rfloor-1$. Then $r\ge L^2/2-2$, and so 
$\left(1+2/L\right)^{r+2}\ge |G|$. From the last statement of (1), $\left(1+2/L\right)^2< 1/2$, and so 
$\left(1+2/L\right)^{r}\ge |G|/2$. From (8), it follows that $|N_i^r|>|G|/2$, and similarly $|N_j^r|>|G|/2$, and so $N_i^r\cap N_j^r\ne \emptyset$. This proves (9).

\bigskip

Let $[t]^{(2)}$ be the set of all two-element
subsets of $\{1\LL t\}$.
Choose $I\subseteq [t]^{(2)}$ maximal such that for each $\{i,j\}\in I$, there is a path $P_{ij}=P_{ji}$ of $G$ with the following 
properties:
\begin{itemize}
\item for each $\{i,j\}\in I$, $P_{ij}$ has one end in $B_i$ and the other in $B_j$, and has no vertices in $X_2$ and no internal vertices in $X_1$;
\item each $P_{i,j}$ is induced and has length at most $L^2-2$; and
\item for all distinct $\{i,j\},\{i',j'\}\in I$, the paths $P_{ij}$ and $P_{i'j'}$ are vertex-disjoint and there are no edges 
between their vertex sets.
\end{itemize}
Let $Y$ be the union of the vertex sets of all the paths $P_{ij}$; then $Y$ is an obstruction set. From (9) and the maximality of $I$,
it follows that $I=[t]^{(2)}$, and so $G$ has an induced subgraph that is a subdivision of $K_t$ with lengths at least three and 
at most $(\log |G|)^2$, a contradiction.
This proves \ref{subdivision}.~\bbox

Since $\Delta(G)\le |G|$, we may replace $d$ by $|G|$, and deduce:
\begin{thm} \label{subdivision2}
Let $G$ be a graph with $|G|\ge 2$, let $k\ge 1$ be an integer with $\omega(G)\le k$, and let $t\ge 3$ be an integer
such that no induced subgraph of $G$
is a subdivision of $K_t$ with lengths between $3$ and $(\log |G|)^2$.  Then
$$\alpha(G)\ge \frac{|G|}{(2t)^{4(k-1)}(\log |G|)^{3k-2}}.$$
\end{thm}
Since a subdivision of $K_t$ with lengths at least two contains a subdivision of any $t$-vertex graph, we could replace $K_t$ by a general graph $H$, and deduce a strengthened version of \ref{subdivision4again}:
\begin{thm} \label{subdivision4}
Let $G$ be a graph with $|G|\ge 2$, let $k\ge 1$ be an integer with $\omega(G)\le k$, and let $H$ be a graph with $t$ vertices, 
such that no induced subgraph of $G$
is a subdivision of $H$ with lengths between $3$ and $(\log |G|)^2$.  Then
$$\alpha(G)\ge \frac{|G|}{(2t)^{4(k-1)}(\log |G|)^{3k-2}}.$$
\end{thm}
If we drop the lower bound on the lengths of the subdivision, then this automatically excludes large cliques, so we can omit the condition that $\omega(G)\le k$. We obtain
a strengthened form of \ref{stable2}:
\begin{thm} \label{subdivision3}
Let $G$ be a graph with $|G|\ge 2$, and let $t\ge 3$ be an integer
such that no induced subgraph of $G$
is a subdivision of $K_t$ with length at most  $(\log |G|)^2$.  Then
$$\alpha(G)\ge \frac{|G|}{(2t)^{4(t-2)}(\log |G|)^{3t-5}}.$$
\end{thm}

\section{Excluding trees and excluding subdivisions}

So we have a pair of theorems: for every forest $H$, all $H$-free graphs with bounded clique number have a nearly-linear stable set 
(proved in~\cite{stable1}), and for every graph $H$, all $H$-subdivision-free graphs with bounded clique number have a nearly-linear stable set (proved in this paper).
There have been previous examples of such pairings between a theorem about excluding a forest, and a theorem about
excluding subdivisions. Here is one, between results concerning the ``strong Erd\H{o}s-Hajnal property'':
\begin{thm}\label{purepairs1}
{\bf\cite{pure1}} For every forest $H$ there exists $\vare>0$ such that in every $H$-free graph $G$ with $|G|\ge 2$ and maximum degree at most $\vare|G|$,
there exist disjoint subsets $A,B$ of $V(G)$ with $|A|,|B|\ge \vare|G|$ such that there are no edges between $A,B$.
\end{thm}
\begin{thm}\label{purepairs2}
{\bf\cite{pure2}} For every graph $H$ there exists $\vare>0$ such that in every $H$-subdivision-free graph $G$ with $|G|\ge 2$ and maximum degree at most $\vare|G|$,
there exist disjoint subsets $A,B$ of $V(G)$ with $|A|,|B|\ge \vare|G|$ such that there are no edges between $A,B$.
\end{thm}

Let $\tau(G)$ denote the maximum $t$ such that $G$ contains a subgraph isomorphic to $K_{t,t}$.
Here is a second pairing:
\begin{thm}\label{kiersteadpenrice}
{\bf\cite{kiersteadpenrice}} For every forest $H$, there is a function $f$
such that every $H$-free graph $G$ has degeneracy at most $f(\tau(G))$.
\end{thm}
\begin{thm}\label{kuhnosthus}
{\bf\cite{kuhnosthus}} For every graph $H$, there is a function $f$
such that every $H$-subdivision-free graph $G$ has degeneracy at most $f(\tau(G))$.
\end{thm}
And a third pairing, between polynomial versions of the two previous results:
\begin{thm}\label{poly1}
{\bf\cite{poly1, hunter}} For every forest $H$, there is a polynomial $f$
such that every $H$-free graph $G$ has degeneracy at most $f(\tau(G))$.
\end{thm}
\begin{thm}\label{girao}
{\bf\cite{cook, girao,hunter}} For every graph $H$, there is a polynomial $f$
such that every $H$-subdivision-free graph $G$ has degeneracy at most $f(\tau(G))$.
\end{thm}
A fourth pairing is between the Gy\'arf\'as-Sumner conjecture \ref{GSconj} and Scott's false conjecture 
\ref{scottconj} (somewhat worrying for those who believe in the predictive value of such pairings and in the truth of 
the Gy\'arf\'as-Sumner conjecture.)

In view of this, one would ask whether any of these four pairings can be unified. Let $H$ be a graph and let $T$ be a forest of $H$. 
We say that a graph $G$ is {\em $H(T)$-subdivision-free} if no induced subgraph of $G$ is isomorphic to a subdivision of $H$
in which the edges of $T$ are not subdivided. This provides a common weakening of the two previous
hypotheses (being $T$-free and being $H$-subdivision-free), and one might hope to unify the pairs on these lines.
For three of the four pairings above this is not known. (In~\cite{pure2} a version of \ref{purepairs2} was proved for $H(T)$-subdivision-free graphs, but only
when $T$ is a path of $H$, not a general forest.)
But we have been able to prove a (numerically somewhat weaker) unification of \ref{stable2} and \ref{stable1}:
\begin{thm}\label{stable1+2}
For every graph $H$, and every tree $T$ of $H$ with radius at most $r\ge 0$, and every $k\ge 2$, let $q=(r+1)(k-1)$; then
there exists $b>0$ such that
every $H(T)$-subdivision-free graph $G$ with clique number at most $k$ and maximum degree at most $d$ satisfies
$$\alpha(G)\ge \frac{|G|}{d^{b(\log d)^{-\frac1q}}(\log |G|)^{b}}.$$
\end{thm}
The proof contains no new ideas, and is a (rather messy) combination of the proofs of \ref{stable2} and \ref{stable1},  
so we will skip it. It will appear in detail in~\cite{tungthesis}.

Finally, there is an intriguing conjecture of Du and McCarty~\cite{mccarty}. 
\begin{thm}\label{mccartyconj}
{\bf Conjecture: }For every class $\mathcal{H}$ of graphs closed under taking induced subgraphs,
if there is a function $f$ such that every $G\in \mathcal{H}$ has degeneracy at most $f(\tau(G))$,  then
$\chi(G)\le |G|^{o(1)}$ for every triangle-free $G\in \mathcal{H}$.
\end{thm}
Our theorems \ref{stable1} and \ref{stable2} both imply special cases of this conjecture, since the existence of
functions $f$ as in \ref{mccartyconj} is proved in~\cite{poly1,kuhnosthus} respectively.

\bigskip \noindent{\bf Acknowledgement:} We would like to thank the referees for their helpful suggestions.

\end{document}